\documentclass[12pt]{article}

\usepackage{amsmath,amssymb,amsbsy,amsfonts,amsthm,latexsym,
            amsopn,amstext,amsxtra,amscd,stmaryrd,color}
\usepackage{dsfont,mathrsfs}
\usepackage{url}
\usepackage{enumerate}

\usepackage[mathscr]{eucal}

\begin{document}

\newtheorem{theorem}{Theorem}
\newtheorem{lemma}{Lemma}
\newtheorem{corollary}{Corollary}
\newtheorem*{computation}{Computation}
\newtheorem{proposition}{Proposition}

\theoremstyle{definition}
\newtheorem*{definition}{Definition}
\newtheorem*{remark}{Remark}
\newtheorem*{example}{Example}


\def\cA{\mathcal A}
\def\cB{\mathcal B}
\def\cC{\mathcal C}
\def\cD{\mathcal D}
\def\cE{\mathcal E}
\def\cF{\mathcal F}
\def\cG{\mathcal G}
\def\cH{\mathcal H}
\def\cI{\mathcal I}
\def\cJ{\mathcal J}
\def\cK{\mathcal K}
\def\cL{\mathcal L}
\def\cM{\mathcal M}
\def\cN{\mathcal N}
\def\cO{\mathcal O}
\def\cP{\mathcal P}
\def\cQ{\mathcal Q}
\def\cR{\mathcal R}
\def\cS{\mathcal S}
\def\cU{\mathcal U}
\def\cT{\mathcal T}
\def\cV{\mathcal V}
\def\cW{\mathcal W}
\def\cX{\mathcal X}
\def\cY{\mathcal Y}
\def\cZ{\mathcal Z}


\def\sA{\mathscr A}
\def\sB{\mathscr B}
\def\sC{\mathscr C}
\def\sD{\mathscr D}
\def\sE{\mathscr E}
\def\sF{\mathscr F}
\def\sG{\mathscr G}
\def\sH{\mathscr H}
\def\sI{\mathscr I}
\def\sJ{\mathscr J}
\def\sK{\mathscr K}
\def\sL{\mathscr L}
\def\sM{\mathscr M}
\def\sN{\mathscr N}
\def\sO{\mathscr O}
\def\sP{\mathscr P}
\def\sQ{\mathscr Q}
\def\sR{\mathscr R}
\def\sS{\mathscr S}
\def\sU{\mathscr U}
\def\sT{\mathscr T}
\def\sV{\mathscr V}
\def\sW{\mathscr W}
\def\sX{\mathscr X}
\def\sY{\mathscr Y}
\def\sZ{\mathscr Z}


\def\fA{\mathfrak A}
\def\fB{\mathfrak B}
\def\fC{\mathfrak C}
\def\fD{\mathfrak D}
\def\fE{\mathfrak E}
\def\fF{\mathfrak F}
\def\fG{\mathfrak G}
\def\fH{\mathfrak H}
\def\fI{\mathfrak I}
\def\fJ{\mathfrak J}
\def\fK{\mathfrak K}
\def\fL{\mathfrak L}
\def\fM{\mathfrak M}
\def\fN{\mathfrak N}
\def\fO{\mathfrak O}
\def\fP{\mathfrak P}
\def\fQ{\mathfrak Q}
\def\fR{\mathfrak R}
\def\fS{\mathfrak S}
\def\fU{\mathfrak U}
\def\fT{\mathfrak T}
\def\fV{\mathfrak V}
\def\fW{\mathfrak W}
\def\fX{\mathfrak X}
\def\fY{\mathfrak Y}
\def\fZ{\mathfrak Z}


\def\C{{\mathds C}}
\def\F{{\mathds F}}
\def\K{{\mathds K}}
\def\L{{\mathds L}}
\def\N{{\mathds N}}
\def\Q{{\mathds Q}}
\def\R{{\mathds R}}
\def\Z{{\mathds Z}}


\def\eps{\varepsilon}
\def\mand{\qquad\mbox{and}\qquad}
\def\\{\cr}
\def\({\left(}
\def\){\right)}
\def\[{\left[}
\def\]{\right]}
\def\<{\langle}
\def\>{\rangle}
\def\fl#1{\left\lfloor#1\right\rfloor}
\def\rf#1{\left\lceil#1\right\rceil}
\def\le{\leqslant}
\def\ge{\geqslant}
\def\ds{\displaystyle}

\def\xxx{\vskip5pt\hrule\vskip5pt}
\def\yyy{\vskip5pt\hrule\vskip2pt\hrule\vskip5pt}
\def\imhere{ \xxx\centerline{\sc I'm here}\xxx }

\newcommand{\commB}[1]{\marginpar{
\vskip-\baselineskip \raggedright\footnotesize
\itshape\hrule\smallskip
\begin{color}{blue} #1\end{color}
\par\smallskip\hrule}}

\newcommand{\commV}[1]{\marginpar{
\vskip-\baselineskip \raggedright\footnotesize
\itshape\hrule\smallskip
\begin{color}{red} #1\end{color}
\par\smallskip\hrule}}


\def\e{{\bf e}}
\def\fp{{\mathfrak p}}

\def\art#1#2{ ({#2},#1)} 
\def\adisc#1{\sD_{#1}} 


\baselineskip17pt

\title{\sc Piatetski-Shapiro Primes \break in a Beatty Sequence}

\author{
{\sc Victor Z. ~Guo} \\
{Department of Mathematics} \\
{University of Missouri} \\
{Columbia, MO 65211 USA} \\
{\tt zgbmf@math.missouri.edu}}

\date{\empty}

\maketitle

\begin{abstract}
Let $\alpha,\beta$ be real numbers such that $\alpha>1$ is irrational and of finite type, and let $c$ be a real number in the range $1<c<\frac{14}{13}$. In this paper, it is shown that there are infinitely many Piatetski-Shapiro primes $p = \fl{n^c}$ in the non-homogenous Beatty sequence $\big(\fl{\alpha m+\beta}\big)_{m=1}^\infty$. 
\end{abstract}





\newpage

\section{Introduction}
For fixed real numbers $\alpha,\beta$ the associated \emph{non-homogeneous Beatty sequence} is the sequence of integers defined by
$$
\cB_{\alpha,\beta}=\(\fl{\alpha n+\beta}\)_{n=1}^\infty,
$$
where $\fl{t}$ denotes the integer part of any $t\in\R$. Such sequences are also called \textit{generalized arithmetic progressions}. It is known that there are infinitely many prime numbers in the Beatty sequence if $\alpha > 0$ (see, for example, the proof of Ribenboim ~\cite[p.~289]{Ribenboim}). 
Moreover, if $\alpha\ge 1$, then the counting function 
$$
\pi_{\alpha,\beta}(x)= \# \big\{\text{\rm prime~}p \le x : p\in\cB_{\alpha,\beta}\big\}
$$
satisfies the asymptotic relation
$$
\pi_{\alpha,\beta} (x) \sim \frac{x}{\alpha \log x} \qquad \text{~\rm as } x \to \infty. 
$$

The \emph{Piatetski-Shapiro sequences} are sequences of the form
$$
\cN^{(c)}=(\fl{n^c})_{n=1}^\infty\qquad(c>1,~c \not\in\N).
$$
Such sequences have been named in honor of Piatetski-Shapiro, who proved~\cite{PS} that $\cN^{(c)}$ contains infinitely many primes if $c\in(1,\frac{12}{11})$. More precisely, for such $c$ he showed that the counting function
$$
\pi^{(c)}(x)= \# \big\{\text{\rm prime~}p\le x : p\in\cN^{(c)}\big\}
$$
satisfies the asymptotic relation
$$
\pi^{(c)}(x) \sim \frac{x^{1/c}}{c \log x} \qquad \text{~\rm as } x \to \infty.
$$
The admissible range for $c$ in this asymptotic formula has been extended many times over the years and is currently known to hold for all $c\in(1,\frac{243}{205})$ thanks to Rivat and Wu~\cite{RivatWu}. The same result is expected to hold for all larger values of $c$.  We remark that if $c\in(0,1)$ then $\cN^{(c)}$ contains all natural numbers, hence all primes in particular.

Since both sequences $\cB_{\alpha,\beta}$ and $\cN^{(c)}$ contain infinitely many primes in the cases described above, it is natural to ask whether infinitely many primes lie in the the intersection $\cB_{\alpha,\beta}\cap\cN^{(c)}$ in some instances.  In this paper we answer this question in the affirmative for certain values of the parameters $\alpha,\beta,c$. Our main result is the following quantitative theorem.
 
\begin{theorem}
\label{thm:main}
Let $\alpha,\beta \in\R$, and suppose that $\alpha > 1$ is irrational and of finite type. Let $c \in (1, \frac{14}{13})$.
There are infinitely many primes in both the Beatty sequence $\cB_{\alpha,\beta}$ and the Piatetski-Shapiro sequence $\sN^{(c)}$. Moreover, the counting function
$$
\pi_{\alpha,\beta}^{(c)}(x)=\big\{\text{\rm prime~}p\le x:p\in\cB_{\alpha,\beta}\cap\cN^{(c)}\big\}
$$
satisfies
$$
\pi_{\alpha,\beta}^{(c)}(x) = \frac{x^{1/c}}{\alpha c \log x}
+O\(\frac{x^{1/c}}{\log^2x}\),
$$
where the implied constant depends only on $\alpha$ and $c$.
\end{theorem}

\medskip\noindent{\bf Remarks.}
We recall that the type $\tau = \tau(\alpha)$ of the irrational number $\alpha$ is defined by
$$
\tau=\sup\big\{t\in\R:\liminf\limits_{n\to\infty} ~n^t\,\llbracket\alpha n\rrbracket=0\big\},
$$
where $\llbracket t\rrbracket$ denotes the distance from a real number $t$ to the nearest integer. For technical reasons we assume that $\alpha$ is of finite type in the statement of the theorem; however, we expect the result holds without this restriction.  

If $\alpha$ is a rational number, then the Beatty sequence $\cB_{\alpha,\beta}$ is a finite union of arithmetic progressions. In the case, Theorem~\ref{thm:main} also holds (in a wider range of $c$) thanks to the work of Leitmann and Wolke~\cite{LW}, who showed that for any coprime integers $a,d$ with $1\le a\le d$ and any real number $c\in(1,\frac{12}{11})$ the counting function
$$
\pi^{(c)}(x; d, a) = \#\big\{p \le x: p \in \sN^{(c)}\text{~and~}p \equiv a \bmod{d}\big\},
$$
satisfies 
\begin{equation}
\label{eq:asymp}
\pi_c(x; d, a)\sim\frac{x^{1/c}}{\phi(d) \log (x)} \qquad \text{~\rm as } x \to \infty,
\end{equation}
where $\phi$ is the Euler function (a more explicit relation than \eqref{eq:asymp} holds in the shorter range
$1<c<\frac{18}{17}$; see Baker \emph{et al} \cite[Theorem~8]{BBBSW}).

We also remark that our theorem is only stated for real numbers $\alpha>1$, for if $\alpha\in(0,1]$ then the set $\cB_{\alpha,\beta}$ contains all but finitely many natural numbers.

\section{Preliminaries}
\subsection{Notation}

We denote by $\fl{t}$ and $\{t\}$ the integer part and the fractional part of $t$, respectively.
As is customary, we put
$$
\e(t)=e^{2\pi it}\mand
\{ t \}=t-\fl{t}\qquad(t\in\R).
$$
Throughout the paper, we make considerable use of the sawtooth function defined by
$$
\psi(t)=t-\fl{t}-\tfrac12=\{t\}-\tfrac12\qquad(t\in\R)
$$

For the Beatty sequence $\cB_{\alpha,\beta}=(\fl{\alpha n + \beta})_{n=1}^\infty$ 
we systematically denote $a=\alpha^{-1}$ and $b=\alpha^{-1}(1-\beta)$. 
For the Piatetski-Shapiro sequence $(\fl{n^c})_{n=1}^\infty$ we always put
$\gamma=1/c$.

Throughout, the letter $p$ always denotes a prime. 

Implied constants in the symbols $O$ and $\ll$ may depend on the parameters $c$ and $A$ (where obvious) but are absolute otherwise.
We use notation of the form $m\sim M$ as an abbreviation for $M< m\le 2M$.

For any set $E$ of real numbers, we denote by $\cX_E$ the characteristic function of $E$; that is,
$$
\cX_E (n) = 
\begin{cases}
1 &\quad\text{if }n \in E, \\
0 &\quad\text{if }n \not\in E. 
\end{cases}
$$

\subsection{Discrepancy}
The \emph{discrepancy} $D(M)$ of a sequence of (not necessarily distinct) real numbers $a_1,a_2,\ldots,a_M\in[0,1)$ is defined by
\begin{equation}
\label{eq:descr_defn}
D(M)=\sup_{\cI\subseteq[0,1)}\bigg|\frac{V(\cI,M)}{M}-|\cI|\,\bigg|,
\end{equation}
where the supremum is taken over all intervals $\cI$ contained in $[0,1)$, $V(\cI,M)$ is the number of positive integers $m\le M$ such that $a_m\in\cI$, and $|\cI|$ is the length of the interval $\cI$.

For any irrational number $\theta$ the sequence of fractional parts $(\{n\theta\})_{n=1}^\infty$ is uniformly distributed over $[0,1)$ (see, e.g., ~\cite[Example~2.1, Chapter~1]{KuNi}). In the special case that $\theta$ is of finite type, the following more precise statement holds 
(see~\cite[Theorem~3.2, Chapter~2]{KuNi}).

\begin{lemma}
\label{lem:discr_with_type}
Let $\theta$ be a fixed irrational number of finite type $\tau$.  Then, for every $\theta\in\R$ the discrepancy $D_{\theta,\mu}(M)$ of the sequence $(\{\theta m+\mu\})_{m=1}^M$ satisfies the bound
$$
D_{\theta,\mu}(M)\le M^{-1/\tau+o(1)}\qquad(M\to\infty),
$$
where the function implied by $o(\cdot)$ depends only on $\theta$.
\end{lemma}

\subsection{Lemmas}

The following lemma provides a convenient characterization of the numbers that occur in the Beatty sequence $\cB_{\alpha,\beta}$.

\begin{lemma}
\label{lem:Beatty_values} 
Let $\alpha,\beta\in\R$ with $\alpha>1$. Then
$$
n\in\cB_{\alpha,\beta}\quad\Longleftrightarrow\quad\cX_a(an+b)=1
$$
where $\cX_a$ is the periodic function defined by
$$
\cX_a(t)=\cX_{(0,a]}(\{t\}) =\begin{cases}
1& \quad \hbox{if $0<\{t\}\le a$}, \\
0& \quad \hbox{otherwise}.
\end{cases}
$$
\end{lemma}

By a classical result of Vinogradov (see \cite[Chapter~I, Lemma~12]{Vin}) we have the following approximation of $\cX_a$ by a Fourier series.
\begin{lemma}
\label{lem:Vinogradov}
For any $\Delta\in(0,\tfrac18)$ with $\Delta\le\tfrac12\min\{a,1-a\}$, there is a real-valued function $\Psi$ with the following properties:
\begin{itemize}
\item[$(i)$~~] $\Psi$ is periodic with period one;
\item[$(ii)$~~] $0 \le\Psi(t)\le 1$ for all $t\in\R$;
\item[$(iii)$~~] $\Psi(t)= \cX_a(t)$ if $\Delta\le \{t\}\le a-\Delta$ or if $a+\Delta\le \{t\}\le 1-\Delta$;
\item[$(iv)$~~] $\Psi(t)=\sum_{k\in\Z}g(k)e(kt)$ for all $t\in\R$,
where $g(0)=a$, and the other Fourier coefficients satisfy the uniform bound
\begin{equation}
\label{eq:coeffbounds}
g(k)\ll\min\big\{|k|^{-1},|k|^{-2}\Delta^{-1}\big\}
\qquad(k\ne 0).
\end{equation}
\end{itemize}
\end{lemma}

We need the following well known approximation of Vaaler~\cite{Vaal}.

\begin{lemma}
\label{lem:Vaaler}
For any $H\ge 1$ there are numbers $a_h,b_h$ such that
$$
\bigg|\psi(t)-\sum_{0<|h|\le H}a_h\,\e(th)\bigg|
\le\sum_{|h|\le H}b_h\,\e(th),\qquad
a_h\ll\frac{1}{|h|}\,,\qquad b_h\ll\frac{1}{H}\,.
$$
\end{lemma}

Next, we recall the following identity for the von Mangoldt function $\Lambda$, which is due to Vaughan (see Davenport~\cite[p.~139]{Daven}).

\begin{lemma}
\label{lem:vaughan}
Let $U,V\ge 1$ be real parameters.  For any $n>U$ we have
$$
\Lambda(n)=
-\sum_{k\,|\,n}a(k)
+\sum_{\substack{cd=n\\d\le V}}(\log c)\mu(d)
-\sum_{\substack{kc=n\\k>1\\c>U}}\Lambda(c)b(k),
$$
where
$$
a(k)=\sum_{\substack{cd=k\\c\le U\\d\le V}}\Lambda(c)\mu(d)
\mand b(k)=\sum_{\substack{d\,\mid\,k\\d\le V}}\mu(d)
$$
\end{lemma}

We also need the following standard result; see~\cite[p.~48]{GraKol}.
\begin{lemma}
\label{lem:index}
For a bounded function $g$ and $N'\sim N$ we have
$$
\sum_{N < p\le N'} g(p) \ll \frac{1}{\log N} \max_{N_1 \le 2N}
\bigg|\sum_{N < n \le N_1} \Lambda(n)g(n)\bigg| + N^{1/2}. 
$$
\end{lemma}

We use the following result of Banks and Shparlinski~\cite[Theorem ~4.1]{BaSh}.

\begin{lemma}
\label{lem:non-liou-theta} Let $\theta$ be a fixed irrational
number of finite type $\tau<\infty$.  Then, for every real number
$0<\eps<1/(8\tau)$, there is a number $\eta>0$ such that the bound
$$
\bigg|\sum_{m\le M}\Lambda(qm+a)\,\e(\theta k m)\bigg|\le
M^{1-\eta}
$$
holds for all integers $1\le k\le M^\eps$ and $0\le a<q\le
M^{\eps/4}$ with $\gcd(a,q)=1$ provided that $M$ is sufficiently
large.
\end{lemma}

We need the following lemma by Van der Corput; see~\cite[Theorem ~2.2]{GraKol}.

\begin{lemma}
\label{lem:GK3}
Let $f$ be three times continuously differentiable on a subinterval
$\cI$ of $(N,2N]$. Suppose that for some $\lambda>0$, the inequalities
$$
\lambda\ll|f''(t)|\ll\lambda\qquad(t\in\cI)
$$
hold, where the implied constants are independent
of $f$ and $\lambda$. Then
$$
\sum_{n\in\cI}\e(f(n))\ll N\lambda^{1/2}+\lambda^{-1/2}.
$$
\end{lemma}

We also need the following two lemmas for the bounds of certain type I and II sums. The two lemmas can be derived by revising the last three lines from the proofs of Baker \emph{et al} \cite[Lemma~24]{BBBSW} and \cite[Lemma~25]{BBBSW}, optimizing the ranges of $K$ and $L$. Specifically we replace $1/3$ and $2/3$ into $3/7$ and $4/7$, respectively. 
\begin{lemma}
\label{sumI}
Suppose $| a_k | \le 1$ for all $k \sim K$. Fix $\gamma \in (0,1)$ and $m, h, d \in \mathbb{N}$. Then for any $K \ll N^{3/7}$ the type I sum
$$
S_I =  \mathop{\sum_{k \sim K} \sum_{l\sim L}}_{ N < kl \le N_1} a_k \: \e(mk^\gamma l^\gamma + klh/d)
$$
satisfies the bound
$$
S_I \ll m^{1/2} N^{3/7 + \gamma/2} + m^{-1/2} N^{1-\gamma/2}.
$$
\end{lemma}

\begin{lemma}
\label{sumII}
Suppose $| a_k | \le 1$ and $| b_l | \le 1$ for $(k,l) \sim (K, L)$. Fix $\gamma \in (0,1)$ and $m, h, d \in \mathbb{N}$. For any $K$ in the range $N^{3/7} \ll K \ll N^{1/2}$, the type II sum
$$
S_{I\!I} = \mathop{\sum_{k \sim K} \sum_{l\sim L}}_{ N < kl \le N_1} a_k b_l \: \e(mk^\gamma l^\gamma + klh/d)
$$
satisfies the bound
$$
S_{I\!I} \ll m^{-1/4}N^{1-\gamma/4} + m^{1/6}N^{16/21+\gamma/6} + N^{11/14}.
$$
\end{lemma}

Finally, we use the following lemma, which provides a characterization of the numbers that
occur in the Piatetski-Shapiro sequence $\cN^{(c)}$.

\begin{lemma}
\label{lem:PS}
A natural number $m$ has the form $\fl{n^c}$ if and only if $\cX^{(c)}(m) = 1$, where
$\cX^{(c)}(m)= \fl{-m^\gamma} - \fl{-(m+1)^\gamma}$.  Moreover,
$$
\cX^{(c)}(m)=\gamma m^{\gamma-1}+\psi(-m^\gamma)-\psi(-(m+1)^\gamma)+O(m^{\gamma-2}).
$$
\end{lemma}

In particular, for any $c\in(1,\frac{243}{205})$ the results of \cite{RivatWu}
yield the estimate
\begin{equation}
\label{eq:PSthm}
\pi^{(c)}(x) = \sum_{p \le x} \cX^{(c)}(p) = \frac{x^\gamma}{c \log x}
+ O\bigg(\frac{x^\gamma}{\log^2 x}\bigg).
\end{equation}

\section{Construction}

In what follows, we use $\tau$ to denote the (finite) type of $\alpha$.

To begin, we express $\pi_{\alpha,\beta}^{(c)}(x)$ as a sum with the characteristic functions of the Beatty and Piatetski-Shapiro sequences; using 
Lemmas~\ref{lem:Beatty_values} and~\ref{lem:PS} we have
$$
\pi_{\alpha,\beta}^{(c)} (x) = \sum_{p \le x}  \cX_a (ap+b) \cX^{(c)}(p).
$$
In view of the properties $(i)$--$(iii)$ of Lemma~\ref{lem:Vinogradov} it follows that
\begin{equation}
\label{eq:basic_estimate}
\pi_{\alpha,\beta}^{(c)} (x)= \sum_{p  \le x} \Psi(ap+b) \cX^{(c)}(p) + O(V(\cI,x))
\end{equation}
holds with some small $\Delta>0$, where $V(\cI,x)$ is the number of primes
$p\in\cN^{(c)}$ not exceeding $x$ for which
$$
\{ap+b\} \in \cI= [0,\Delta) \cup (\alpha-\Delta,\alpha+\Delta) \cup (1-\Delta,1);
$$
that is,
$$
V(\cI,x) = \sum_{p \le x} \cX_{\cI}\big(\{ap+b\}\big) \cX^{(c)}(p).
$$
By Lemma~\ref{lem:PS} we see that
$$
V(\cI,x)=\gamma V_1(x)+V_2(x)+O(1),
$$
where
\begin{align*}
V_1(x) &= \sum_{p \le x} \cX_{\cI}\big(\{ap+b\}\big) p^{\gamma-1},\\
V_2(x) &= \sum_{p \le x} \cX_{\cI}\big(\{ap+b\}\big) \big(\psi(-p^\gamma) - \psi(-(p+1)^\gamma)\big).
\end{align*}
Using \eqref{eq:PSthm} we immediately derive the bound
\begin{align*}
V_2(x) &\le \sum_{p \le x} \big(\psi(-p^\gamma) - \psi(-(p+1)^\gamma)\big) \ll \frac{x^\gamma}{\log^2 x}\,.
\end{align*}
To bound $V_1(x)$ we split the sum over $n \le x$ into $O(\log x)$ dyadic intervals of the form $(N, 2N]$ with $N \ll x$ and apply Lemma~\ref{lem:index}, obtaining that
\begin{align*}
V_1(x) 
&\ll \log x\cdot \max_{N\le x}\bigg(\frac{1}{\log N} \max_{N_1 \le 2N}
\bigg|\sum_{N < n \le N_1} \Lambda(n) X_{\cI} (\{an+b\})
n^{\gamma-1}\bigg| + N^{1/2} \bigg)  \\
&\ll x^{\gamma-1} \log x\cdot \max_{N\le x}\max_{N_1 \le 2N} \bigg|
\sum_{N < n < N_1} X_{\cI} (\{an+b\})\bigg|+x^{1/2}\log x.
\end{align*}
Since $|\cI|=4\Delta$, it follows from the definition \eqref{eq:descr_defn} and
Lemma~\ref{lem:discr_with_type} that
$$
V_1(x) \ll \Delta x^\gamma \log x +x^{\gamma-\frac{1}{\tau}+o(1)}\qquad(x\to\infty).
$$
Therefore, 
\begin{equation}
\label{eq:bound V(I,x)} V(\cI,x) \ll \Delta x^\gamma \log x+\frac{x^\gamma}{\log^2 x}\,.
\end{equation}

Now let $K\ge\Delta^{-1}$ be a large real number, and let $\Psi_K$ be the trigonometric polynomial defined by
\begin{equation}
\label{eq:PKdefn} \Psi_K(t)=\sum_{|k|\le K}g(k)e(kt).
\end{equation}
Using \eqref{eq:coeffbounds} it is clear that the estimate
\begin{equation}
\label{eq:PKP} \Psi(t)=\Psi_K(t)+O(K^{-1}\Delta^{-1})
\end{equation}
holds uniformly for all $t\in\R$. Combining \eqref{eq:PKP} with \eqref{eq:basic_estimate} and taking into
account \eqref{eq:bound V(I,x)} we derive that
$$
\pi_{\alpha,\beta}^{(c)} (x)= \sum_{p \le x} \Psi_K(ap+b) \cX^{(c)}(p)+ O(E(x)),
$$
where 
$$
E(x) = \Delta x^\gamma \log x+\frac{x^\gamma}{\log^2 x}
+ K^{-1} \Delta^{-1} \sum_{p  \le x}  \cX^{(c)}(p).
$$
For fixed $A\in(0,1)$ we put
$$
\Delta=x^{-A/2}\mand K=x^A.
$$
Note that our previous application of Lemma~\ref{lem:Vinogradov} to deduce
\eqref{eq:basic_estimate} is justified.  Use these values of
$\Delta$ and $K$ along with \eqref{eq:PSthm} we obtain that
$$
E(x) \ll x^{\gamma-A/2} \log x+ \frac{x^{\gamma}}{\log^2 x} +\frac{x^{\gamma- A/2}}{\log x}
\ll \frac{x^{\gamma}}{\log^2 x}\,.
$$
Using the definition \eqref{eq:PKdefn} it therefore follows that
\begin{equation}
\label{eq:facebook}
\pi_{\alpha,\beta}^{(c)} (x) = \sum_{p \le x} \sum_{|k| \le x^A} g(k) \e(kap+kb) \cX^{(c)}(p)+ O\(\frac{x^{\gamma}}{\log^2 x}\).
\end{equation}

Next, using Lemma~\ref{lem:PS} we express the double sum in \eqref{eq:facebook}
as $\sum_1 + \sum_{2,1} + \sum_{2,2}$ with
\begin{align*}
{\textstyle\sum_1} &=  g(0)\sum_{p \le x} \cX^{(c)}(p),\\
{\textstyle\sum_{2,1}} &= \sum_{\substack{k \neq 0 \\ |k| \le x^A}} g(k) \sum_{p \le x} \e(kap+kb) \big( \gamma p^{\gamma-1} + O(p^{\gamma-2}) \big), \\
{\textstyle\sum_{2,2}} &= \sum_{\substack{k \neq 0 \\ |k| \le x^A}} g(k) \sum_{p \le x} \e(kap+kb) \big\{ \psi(-(p+1)^\gamma) - \psi(-p^\gamma) \big\}.
\end{align*}
Recalling that $g(0) = \alpha^{-1}$ we have
$$
{\textstyle\sum_1}= \alpha^{-1} \sum_{p \le x} \cX^{(c)}(p) = \frac{x^\gamma}{\alpha c \log x} + O\(\frac{x^\gamma}{\log^2 x}\),
$$
which provides the main term in our estimation of $\pi_{\alpha,\beta}^{(c)}(x)$. 

To bound ${\textstyle\sum_{2,1}}$ we follow the method used above to bound $V(\cI,x)$ and use partial summation together with \eqref{eq:coeffbounds} to conclude that
$$
{\textstyle\sum_{2,1}}
\ll x^{\gamma-1}\log x\sum_{\substack{k \neq 0 \\ |k| \le x^A}}
\frac{1}{|k|} \max_{N\le x}\bigg(\frac{1}{\log N} \max_{N' \le 2N}
\bigg|\sum_{N \le n \le N'} \Lambda(n) \e(k \alpha^{-1} n) \bigg|+ 1 \bigg) 
$$
Assuming as we may that $0 < A < 1/(8\tau)$, by Lemma~\ref{lem:non-liou-theta} it follows that there exists $\eta \in (0,1)$ such that the bound
$$
\max_{N\le x}\bigg(\frac{1}{\log N} \max_{N' \le 2N}
\bigg|\sum_{N \le n \le N'} \Lambda(n) \e(k \alpha^{-1} n) \bigg|\bigg)\ll x^{1-\eta}
$$
holds uniformly for $|k|\le x^A$, $k\ne 0$.  Consequently, we derive the bound
$$
{\textstyle\sum_{2,1}}
\ll  \big( x^{\gamma-1} x^{1-\eta} + x^{\gamma-1} \big) \log^2 x
\ll \frac{x^\gamma}{\log^2 x}\,,
$$
which is acceptable.

To complete the proof it suffices to show that
${\textstyle\sum_{2,2}}\ll x^\gamma/\log^2 x$.
To accomplish this task we use the method in~\cite[pp.~47--53]{GraKol}.  Denote
$$
{\textstyle\sum_3} = \sum_{p \le x} \e(kap+kb) \big\{ \psi(-(p+1)^\gamma) - \psi(-p^\gamma) \big\}.
$$
It is enough to show that the bound $\sum_3 \ll x^{\gamma - \eps}$ holds with some $\eps>0$ uniformly for $k$, for then we have by \eqref{eq:coeffbounds}:
$$
{\textstyle \sum_{2,2}} \ll \sum_{\substack{k \neq 0 \\ |k| \le x^A}} \frac{1}{|k|} \cdot x^{\gamma - \eps} \ll x^{\gamma - \eps} \log x \ll
\frac{x^\gamma}{\log^2x}\,.
$$
By Lemma ~\ref{lem:Vaaler}, for any $H \ge 1$ we can write
$$
{\textstyle\sum_3} = {\textstyle\sum_4} + O\big({\textstyle\sum_5}\big),
$$
where
\begin{align*}
{\textstyle\sum_4} &= \sum_{p \le x}  \sum_{0 < |h| \le H} a_h \bigl( \e(kap+kb+h(p+1)^\gamma) - \e(kap+kb+hp^\gamma) \bigr),\\
{\textstyle\sum_5} &= \sum_{n \le x} \sum_{|h|\le H}b_h \bigl(\e(kan+kb+h(n+1)^\gamma)+\e(kan+kb+hn^\gamma)\bigl),
\end{align*}
with some numbers $a_h,b_h$ that satisfy $a_h \ll |h|^{-1}$ and $b_h \ll H^{-1}$.
Thus, it suffices to show that the bounds $\sum_4 \ll x^{\gamma-\eps}$ and $\sum_5 \ll x^{\gamma-\eps}$ hold  with an appropriate choice of $H$.
To this end, we put
$$
H=x^{1-\gamma+2\eps}.
$$

First, we consider $\sum_5$. The contribution from $h=0$ is 
\begin{equation}
\label{eq:h=0}
2\sum_{n < x} b_0 \e(kan + kb) \ll b_0 {|ka|}^{-1} \ll 1. 
\end{equation}
Suppose that $N \le x$ and $N_1\sim N$. We denote 
$$
S_j = \sum_{N < n \le N_1} \sum_{0 < |h| \le H}b_h \e(kan+kb+h(n+j)^\gamma).
$$
To bound the part that $h \neq 0$, it is suffices to show that $S_j \ll x^{1-\eps}$ for $j = 0 \text{ or } 1$. By a shift of $n$, we have 
$$
S_j \ll \sum_{N < n \le N_1} H^{-1} \sum_{0 < h \le H} \e(kan+hn^\gamma).
$$ 
Using Lemma~\ref{lem:GK3} with the choice of $\lambda = hN^{\gamma - 2}$, we obtain
\begin{align*}
S_j
&\ll H^{-1} \sum_{0 < h \le H} \big( N (hN^{\gamma - 2})^{1/2} + (hN^{\gamma - 2})^{-1/2} \big) \\
&\ll (x^{1-\gamma+2\eps})^{1/2} x^{\gamma/2} + (x^{1-\gamma+2\eps})^{-1/2} x^{1-\gamma/2} \ll x^{1/2 + 2\eps}.
\end{align*}
Then summing over $N$, adding the part that $h=0$ from~\eqref{eq:h=0} and recalling that $\gamma>1/2$, we see that the bound
$$
{\textstyle\sum_5} \ll x^{1/2+2\eps} \log x + 1 \ll x^{\gamma-\eps}
$$
holds if the parameter $\eps$ is sufficiently small, which we can assume. 

To bound $\sum_4$ we apply Lemma~\ref{lem:index} and split the sum into $O(\log x)$ 
dyadic intervals of $(N, N_1]$ to derive the bound
\begin{align*}
&\sum_{N < p \le N_1}\sum_{0 < |h| \le H} a_h \bigl( \e(kap + kb + h(p+1)^\gamma) - \e(kap+kb+hp^\gamma) \bigr) \\
&\qquad \ll \frac{N^{\gamma-1}}{\log N} \max_{N_2\le 2N} \biggl| \sum_{1\le h\le H} \sum_{N < n \le N_2} \Lambda(n) \e(kan+kb + hn^\gamma)\biggl|+N^{1/2}.
\end{align*}
Summing over $N$ and taking into account that $\gamma>1/2$,
we obtain the desired bound $\sum_4\ll x^{\gamma}/\log^2x$
(hence also $\sum_3\ll x^{\gamma}/\log^2x$) provided that
\begin{equation}
\label{eq:sumlamF}
\sum_{1\le h\le H} \sum_{N<n\le N_2} \Lambda(n)\e(kan+kb + hn^\gamma) \ll x^{1-\eps}.
\end{equation}
Using Lemma~\ref{lem:vaughan}, we can express the sum on the left side of \eqref{eq:sumlamF} as
$$
\sum_{1\le h\le H}(-S_{1,h}+S_{2,h}-S_{3,h}),
$$
where
\begin{align*}
S_{1,h}
&=\sum_{m\le UV} \sum_{N/m\le n\le N_2/m} 
\tilde a(m) \e(kamn + kb + hm^\gamma n^\gamma),\\
S_{2,h}
&=\sum_{m\le V}\sum_{N/m\le n\le N_2/m} 
\mu(m)(\log n)\e(kamn + kb + hm^\gamma n^\gamma),\\
S_{3,h}
&=\sum_{V<n<N_2/U} \sum_{\substack{N/n\le m\le N_2/n\\m>U}}
\tilde b(n) \Lambda(m) \e(kamn +kb + hm^\gamma n^\gamma),
\end{align*}
and the functions $\tilde a$ and $\tilde b$ are given by
$$
\tilde a(m)=\sum_{\substack{cd=m\\c\le U\\d\le V}}\Lambda(c)\mu(d)
\mand \tilde b(n)=\sum_{\substack{d\,\mid\,n\\d\le V}}\mu(d).
$$
To establish \eqref{eq:sumlamF} it suffices to show that
\begin{equation}
\label{eq:sumSjs}
\sum_{1\le h\le H}S_{j,h}\ll x^{1-\eps}\qquad(j=1,2,3).
\end{equation}

We turn to the problem of bounding $S_{1,h}$, $S_{2,h}$ and $S_{3,h}$.
The sum $S_{2,h}$ is of type I, and $S_{3,h}$ is of type II.
To bound $S_{1,h}$ we write it in the form $S_{4,h}+S_{5,h}$, where
$S_{4,h}$ is a type I sum and $S_{5,h}$ is a type II sum.
To simplify the calculation, we take 
$$
V=N^{3/7} \mand U=N^{1/7}.
$$
Since $V \ll N^{3/7}$, we apply Lemma~\ref{sumI} to bound the sum $S_{2,h}$. 
\begin{align*}
\sum_{1 \le h \le H} S_{2,h} 
&\ll \sum_{1 \le h \le H} \log N |\sum_{m \le V} \sum_{N/m \le n \le {N_2}/m} \e(kamn + hm^\gamma n^\gamma) | \\
&\ll \sum_{1 \le h \le H} \log N (h^{1/2} N^{3/7 + \gamma/2} + h^{-1/2} N^{1-\gamma/2} ) \\
&\ll x^{27/14 - \gamma +3\eps} + x^{3/2 - \gamma + \eps} \ll x^{1- \eps}
\end{align*}
if assuming $\gamma > \frac{13}{14}$.

The sum $S_{3,h}$ can be split into $\ll \log^2 N$ subsums of the form
$$
\sum_{X \le m \le 2X} \sum_{\substack{Y \le n \le 2Y \\ N \le mn \le N_1}} \alpha(m) \beta(n) \e (k \alpha^{-1} mn + h m^\gamma n^\gamma).
$$
It suffices to consider the special case that $V < Y \le N^{1/2}$ and $ N^{1/2} < X \le N/V$. Applying Lemma~\ref{sumII}
(taking into account the estimates
$\alpha(m)\ll N^{\eps/2}$ and $\beta(n)\ll N^{\eps/2}$)
each subsum is
$$
\ll \( h^{-1/4} N^{1-\gamma/4} + h^{1/6} N^{16/21 + \gamma/6} + N^{11/14} \) N^{\eps}.
$$
Therefore, the bound
\begin{align*}
\sum_{1\le h\le H}S_{3,h}
&\ll\( H^{3/4} N^{1-\gamma/4} + H^{7/6} N^{16/21 + \gamma/6} + HN^{11/14} \) N^{\eps}\\
&\ll\((x^{1-\gamma+2\eps})^{3/4} x^{1-\gamma/4} + (x^{1-\gamma+2\eps})^{7/6} x^{16/21 + \gamma/6} + (x^{1-\gamma+2\eps})x^{11/14} \) x^{\eps}\\
&\ll\(x^{7/4-\gamma} + x^{27/14-\gamma}+ x^{25/14-\gamma}\) x^{4\eps} \ll x^{1-\eps}
\end{align*}
under our hypothesis that $\gamma > \frac{13}{14}$. 

Finally, to derive the required bound $S_{1,h}\ll x^{1-\eps}$ we write
$$
S_{1,h} = S_{4,h} + S_{5,h},
$$
where
\begin{align*}
S_{4,h} &=\sum_{m \le V} \sum_{N/m \le n \le N_2/m} a(m)
\e(kamn + kb + h m^\gamma n^\gamma), \\
S_{5,h} &=\sum_{V < m \le UV} \sum_{N/m \le n \le N_2/m} a(m)
\e(kamn + kb + h m^\gamma n^\gamma).
\end{align*}
Since $a(m)\le\log m$ the methods used above to bound $S_{2,h}$ and $S_{3,h}$
can be applied to $S_{4,h}$ and $S_{5,h}$, respectively, to see that the bounds
\begin{equation}
\label{eq:sumSjs}
\sum_{1\le h\le H}S_{j,h}\ll x^{1-\eps}\qquad(j=4,5).
\end{equation}
hold under our hypothesis that $\gamma>\frac{13}{14}$.
This establishes \eqref{eq:sumSjs}, and the theorem is proved.

\section{Remarks}

We note that both \cite[Theorem~7]{BBBSW} and \cite[Theorem~8]{BBBSW} can be improved using Lemma~\ref{sumI} and Lemma~\ref{sumII} instead of \cite[Lemma~24]{BBBSW} and \cite[Lemma~25]{BBBSW}, respectively. The range of $c$ in \cite[Theorem~7]{BBBSW} can be extended from $(1,\frac{147}{145})$ to $(1,\frac{571}{561})$, with a small improvement of $0.004$. For \cite[Theorem~8]{BBBSW}, the range of $c$ is improved from $(1,\frac{18}{17})$ to $(1,\frac{14}{13})$ and the error term is improved from $O(x^{17/39 + 7\gamma/13 + \eps})$ to $O(x^{3/7 + 7\gamma/13 + \eps})$. 

It would be interesting to see whether the range of $c$ in the statement of Theorem~\ref{thm:main} can be improved using more sophisticated methods to improve our type II estimates. With more work, it should be possible to remove our assumption that $\alpha$ is of finite type.  For the sake of simplicity, these ideas have not been pursued in the present paper.

\bigskip\noindent{\bf Acknowledgement.}
The author would like to thank his advisor, William Banks, for suggesting this work and for several helpful discussions.


\begin{thebibliography}{99}

\bibitem{BBBSW}
R.~C.~Baker, W.~D.~Banks, J.~Br\"udern, I.~E.~Shparlinski and A.~J.~Weingartner,
`Piatetski-Shapiro sequences,'
\emph{Acta Arith.} \textbf{157} (2013), no.~1, 37--68. 

\bibitem{BaSh} 
W.~D.~Banks and I.~E.~Shparlinski,
`Prime numbers with Beatty sequences,' \emph{Colloq.\ Math.}
\textbf{115}  (2009),  no.~2, 147--157.

\bibitem{Daven}
H.~Davenport
\emph{Multiplicative number theory}. 
Graduate Texts in Mathematics, \textbf{74}. Springer-Verlag, New York-Berlin, 1980.

\bibitem{GraKol}
S.~W.~Graham and G.~Kolesnik,
\emph{Van der Corput's method of exponential sums}.
London Mathematical Society Lecture Note Series, \textbf{126}.
Cambridge University Press, Cambridge, 1991.

\bibitem{KuNi}
L.~Kuipers and H.~Niederreiter, \emph{Uniform distribution of
sequences}. Pure and Applied Mathematics. Wiley-Interscience, New
York-London-Sydney, 1974.

\bibitem{PS}
I.~I.~Piatetski-Shapiro, `On the distribution of prime numbers
in the sequence of the form $\fl{f(n)}$,'
\emph{Mat. Sb.} \textbf{33} (1953), 559--566.

\bibitem{Ribenboim}
P.~Ribenboim
\emph{The new book of prime number records}.
Springer-Verlag, New York, 1996.

\bibitem{RivatWu}
J.~Rivat and J.~Wu,
`Prime numbers of the form $\fl{n^c}$,'
\emph{Glasg.\ Math.\ J.} \textbf{43} (2001), no.~2, 237--254. 

\bibitem{Vaal}
J.~D.~Vaaler,
`Some extremal problems in Fourier analysis,'
\emph{Bull.\ Amer.\ Math.\ Soc.} \textbf{12} (1985), 183--216. 

\bibitem{Vin} I.~Vinogradov,
\emph{The method of trigonometrical sums in the theory of numbers}.
Dover Publications, Inc., Mineola, NY, 2004.

\bibitem{LW}
D.~Leitmann and D.~Wolke,
`Primzahlen der Gestalt $[n^\Gamma]$ in arithmetischen progressionen', (German) 
\emph{Arch. \ Math.} (Basel) \textbf{25} (1974), 492-494.


\end{thebibliography}
\end{document}